# LIMITING SHAPE FOR DIRECTED PERCOLATION MODELS

By James B. Martin

*CNRS and Université Paris 7*

We consider directed first-passage and last-passage percolation on the nonnegative lattice $\mathbb{Z}_+^d$, $d \geq 2$, with i.i.d. weights at the vertices. Under certain moment conditions on the common distribution of the weights, the limits $g(\mathbf{x}) = \lim_{n \to \infty} n^{-1} T(\lfloor n\mathbf{x} \rfloor)$ exist and are constant a.s. for $\mathbf{x} \in \mathbb{R}_+^d$, where $T(\mathbf{z})$ is the passage time from the origin to the vertex $\mathbf{z} \in \mathbb{Z}_+^d$. We show that this *shape function* $g$ is continuous on $\mathbb{R}_+^d$, in particular at the boundaries. In two dimensions, we give more precise asymptotics for the behavior of $g$ near the boundaries; these asymptotics depend on the common weight distribution only through its mean and variance. In addition we discuss growth models which are naturally associated to the percolation processes, giving a shape theorem and illustrating various possible types of behavior with output from simulations.

**1. Introduction.** We consider directed *first-passage* and *last-passage* percolation models on the nonnegative lattice $\mathbb{Z}_+^d$, focusing in particular on behavior close to the boundaries of the orthant.

With each node $\mathbf{z} \in \mathbb{Z}_+^d$, associate the *weight* $X(\mathbf{z})$. We assume that the weights $\{X(\mathbf{z}), \mathbf{z} \in \mathbb{Z}_+^d\}$ are i.i.d. according to some common distribution $F$ on $\mathbb{R}$; by allowing the weights to take negative as well as positive values we can consider first-passage and last-passage models simultaneously.

A *directed path* in $\mathbb{Z}_+^d$ is a path each step of which consists of increasing a single coordinate by 1. Let $T(\mathbf{z})$, the *last-passage time to* $\mathbf{z}$, be the maximum weight of all directed paths from the origin to the point $\mathbf{z}$, where the weight of a path is the sum of the weights of the nodes along the path. (See Section 2 for precise definitions.) Natural objects of study are asymptotic quantities such as the function $g : \mathbb{R}_+^d \mapsto \mathbb{R}$ defined by

$$g(\mathbf{x}) = \sup_{n \in \mathbb{N}} \frac{\mathbb{E} T(\lfloor n\mathbf{x} \rfloor)}{n}.$$

From superadditivity properties, we have that this supremum is in fact a limit, and that $n^{-1}T(\lfloor n\mathbf{x} \rfloor) \to g(\mathbf{x})$ a.s. as $n \to \infty$, for all $\mathbf{x} \in \mathbb{R}^d_+$. We call $g$ the *shape function*, since it determines the limiting shape for the growth model associated to the percolation process; $g(\mathbf{x})$ is also sometimes known as the *time constant* in direction $\mathbf{x}$.

Analogous first-passage percolation models on the undirected lattice are by now very well known—see, for example, Kesten [14, 15] or Durrett [7] for fundamental results. Recently, the directed last-passage model has also received much attention; in particular, the case where $d = 2$ and the weight distribution $F$ is geometric or exponential. First, this is essentially the only nontrivial case (whether directed or undirected, first- or last-passage) where the form of the shape function $g$ above is known; for exponential weights, it was first given by Rost [25]. But much more precise results are now available; in particular Johansson [13] extended methods developed by Baik, Deift and Johansson [2] for the closely related model of the longest increasing subsequence of a random permutation, and showed that, for $\alpha > 0$, the distribution of $n^{-1/3}\{T((n, \lfloor \alpha n \rfloor)) - ng((1, \alpha))\}$ converges as $n \to \infty$ to a nondegenerate limit (the Tracy–Widom distribution, which also arises as the limiting distribution for the size of the largest eigenvalue of a random matrix from the Gaussian unitary ensemble).

Two-dimensional directed last-passage percolation problems with general weight distributions have also been studied in detail in the context of tandem queueing systems; see, for example, [1, 10, 19].

Our first observation is a condition on the weight distribution $F$ under which the shape function $g$ above is finite everywhere. The condition required on the positive tail is that $\int_0^\infty (1 - F(s))^{1/d} \, ds < \infty$; this follows quickly from analogous results for the related model of *greedy lattice animals* (introduced in [5] and [9]; the precise results we use are from [20]). We note that, as in the greedy lattice animals model, there is still a small gap between this sufficient condition and the best currently known necessary condition, which is that $\mathbb{E} X_+^d < \infty$ (see [20] for a discussion).

Our first main result is then that the shape function $g$ is continuous on all of $\mathbb{R}^d_+$, including at the boundaries (in fact, continuity on the interior is immediate from a simple concavity property). We note that the question of continuity for the directed first-passage model was raised by Newman and Piza [21]; for the two-dimensional last-passage model it was resolved (in a queueing theory context) by Glynn and Whitt [10] for distributions with an exponential tail, and then by Baccelli, Borovkov and Mairesse [1] and Martin [19] under weaker moment conditions. Particular tools which we use to prove continuity at the boundaries in any dimension are a truncation which relies on a bound given in [20] for the growth rate in the greedy lattice animals model, and a concentration inequality derived from a result of Talagrand [27].

In two dimensions we then give more precise information about the behavior of $g$ close to the boundary. For a distribution $F$ with mean $\mu$ and variance $\sigma^2$, we prove the asymptotic formula, as $\alpha \to 0$,

$$g((1,\alpha)) = \mu + 2\sigma\sqrt{\alpha} + o(\sqrt{\alpha}).$$

In addition to the tools used in the proof of continuity, we use here a comparison with a variant form of directed percolation analyzed by Seppäläinen [26], and an estimate for the speed of convergence in the central limit theorem from [24], in order to prove a universality property over all $F$ for the asymptotics at the boundary. Comparing the exact formula known for the case where $F$ is the exponential distribution then yields the result. These asymptotics are linked to the *Brownian directed percolation* model—obtained, loosely speaking, by reversing the order of the limits $n \to \infty$ (in the definition of $g$) and $\alpha \to 0$—which has been widely studied recently in various contexts; see, for example, [3, 11, 12, 23]. See also [22] for a survey of the connections between these various directed percolation processes, random matrix theory and noncolliding particle systems.

For various results on the dependence of the time constant on the weight distribution in the context of undirected first-passage percolation, see, for example, [16, 18, 28].

Just as in the case of undirected first-passage percolation, there are *growth processes* naturally associated to the directed percolation models considered here. In Section 5 we describe these and prove a *shape theorem* analogous to those given in [4] and [14]. We also discuss, with illustrations from simulations (see Figures 1–7), various possible behaviors of the growth processes, and the differences which exist between the directed and undirected cases and between the first-passage and last-passage cases.

**2. Notation and main results.** We work with the $d$-dimensional nonnegative lattice $\mathbb{Z}_+^d$. For $\mathbf{x} \in \mathbb{Z}_+^d$ (and similarly $\mathbb{Z}^d$, $\mathbb{R}_+^d$ and $\mathbb{R}^d$) we write $x_i$ for the $i$th component of $\mathbf{x}$; we use the norm $\|\mathbf{x}\| = \sum |x_i|$, and write $\mathbf{x} \leq \mathbf{x}'$ if $x_i \leq x_i'$ for $i = 1, \ldots, d$. We write $\mathbf{0}$ for the origin and $\mathbf{1}$ for the point all of whose coordinates equal 1, and $\mathbf{e}_i$ for the point all of whose coordinates are 0 except the $i$th which is 1 (so that $\mathbf{x} = \sum_{i=1}^d x_i \mathbf{e}_i$).

With each point $\mathbf{v}$ of $\mathbb{Z}_+^d$, associate the *weight* $X(\mathbf{v})$. We assume that the weights $\{X(\mathbf{v}), \mathbf{v} \in \mathbb{Z}_+^d\}$ are i.i.d. random variables, with common distribution function $F$, where $F(s) = \mathbb{P}(X \leq s)$ (as here, we sometimes write simply $X$ to denote a generic random variable with distribution $F$).

A *directed path* in $\mathbb{Z}_+^d$ is a path each step of which consists of increasing a single coordinate by 1. For $\mathbf{z}_1, \mathbf{z}_2 \in \mathbb{Z}_+^d$, with $\mathbf{z}_1 \leq \mathbf{z}_2$, let $\Pi[\mathbf{z}_1, \mathbf{z}_2)$ be the set of directed paths from $\mathbf{z}_1$ to $\mathbf{z}_2$. We identify a path with the set of points it contains; by convention we exclude the final point $\mathbf{z}_2$ (but include the initial point $\mathbf{z}_1$, unless $\mathbf{z}_1 = \mathbf{z}_2$). Note that all paths in $\Pi[\mathbf{z}_1, \mathbf{z}_2)$ have size (or "length") $\|\mathbf{z}_2 - \mathbf{z}_1\|$.

For $\mathbf{z}_1 \leq \mathbf{z}_2$, define $T(\mathbf{z}_1, \mathbf{z}_2)$, the *last-passage time* from $\mathbf{z}_1$ to $\mathbf{z}_2$, by

$$T(\mathbf{z}_1, \mathbf{z}_2) = \max_{\pi \in \Pi[\mathbf{z}_1, \mathbf{z}_2)} \sum_{\mathbf{v} \in \pi} X(\mathbf{v}).$$

In the case $\mathbf{z}_1 = \mathbf{0}$, we write simply $\Pi[\mathbf{z}] = \Pi[\mathbf{0}, \mathbf{z})$ and

$$(2.1) \qquad T(\mathbf{z}) = T(\mathbf{0}, \mathbf{z}) = \max_{\pi \in \Pi[\mathbf{z}]} \sum_{\mathbf{v} \in \pi} X(\mathbf{v}).$$

We immediately have the following superadditivity property: if $\mathbf{z}_1 \leq \mathbf{z}_2 \leq \mathbf{z}_3 \in \mathbb{Z}_+^d$, then

$$(2.2) \qquad T(\mathbf{z}_1, \mathbf{z}_2) + T(\mathbf{z}_2, \mathbf{z}_3) \leq T(\mathbf{z}_1, \mathbf{z}_3).$$

Suppose $\mathbb{E}|X| < \infty$. Then also $\mathbb{E}|T(\mathbf{z})| < \infty$ for all $\mathbf{z} \in \mathbb{Z}_+^d$. For $\mathbf{x} \in \mathbb{R}_+^d$, we now define

$$(2.3) \qquad g(\mathbf{x}) = \sup_{n \in \mathbb{N}} \frac{1}{n} \mathbb{E}\, T(\lfloor n\mathbf{x} \rfloor)$$

(which may be infinite). We sometimes write $g_F(\mathbf{x})$ to emphasize the dependence on the distribution $F$; we also write $g(x_1, \ldots, x_d)$ for $g(\mathbf{x})$ when $\mathbf{x} = (x_1, \ldots, x_d)$ in order to avoid proliferation of brackets.

The following basic properties of the function $g$ are immediate from this definition and from the superadditivity in (2.2), using (a superadditive version of) Kingman's subadditive ergodic theorem:

PROPOSITION 2.1. *Suppose* $\mathbb{E}|X| < \infty$.

(i) *For all* $\mathbf{x} \in \mathbb{R}_+^d$,

$$\frac{1}{n} T(\lfloor n\mathbf{x} \rfloor) \to g(\mathbf{x}) \qquad \text{as } n \to \infty, \text{ a.s. and } [\text{if } |g(\mathbf{x})| < \infty] \text{ in } \mathcal{L}_1.$$

(ii) $g(\alpha \mathbf{x}) = \alpha g(\mathbf{x})$ *for all* $\alpha \geq 0$, $\mathbf{x} \in \mathbb{R}_+^d$.
(iii) $g$ *is invariant under permutations of the coordinates.*
(iv) $g(\mathbf{x}) + g(\mathbf{y}) \leq g(\mathbf{x} + \mathbf{y})$ *for all* $\mathbf{x}, \mathbf{y} \in \mathbb{R}_+^d$.

The following result gives conditions under which the function $g$ is finite. Condition (2.5) is stronger than the condition that $\mathbb{E} X_+^d < \infty$ (which is known to be necessary for the finiteness) but weaker, for example, than the condition $\mathbb{E} X_+^d (\log_+ X)^{d-1+\varepsilon} < \infty$. See [20] for details.

PROPOSITION 2.2. *If*

$$(2.4) \qquad \mathbb{E}|X| < \infty$$

and

(2.5) $$\int_0^\infty (1 - F(s))^{1/d}\, ds < \infty,$$

then $|g(\mathbf{x})| < \infty$ for all $\mathbf{x} \in \mathbb{R}_+^d$.

Our first main result is then:

THEOREM 2.3. *Under conditions* (2.4) *and* (2.5), *the shape function $g$ is continuous on all of $\mathbb{R}_+^d$ (including at the boundaries).*

Proposition 2.2 and Theorem 2.3 are proved in Section 3. In Section 4 we prove the following theorem, which gives more precise asymptotics for $g$ at the boundary in the case $d = 2$.

THEOREM 2.4. *Let $d = 2$. Let the distribution $F$ have mean $\mu$ and variance $\sigma^2$, and satisfy*

(2.6) $$\int_0^\infty (1 - F(s))^{1/2}\, ds < \infty$$

and

(2.7) $$\int_{-\infty}^0 F(s)^{1/2}\, ds < \infty.$$

*Then as $\alpha \downarrow 0$,*

(2.8) $$g(1, \alpha) = \mu + 2\sigma\sqrt{\alpha} + o(\sqrt{\alpha}).$$

Note that the framework effectively includes first-passage as well as last-passage percolation models, since the weights may take negative as well as positive values; replacing *max* by *min* and replacing the weights $X(\mathbf{z})$ by $-X(\mathbf{z})$ would simply change the sign of $T$ and so of $g$. When considering associated growth models in Section 5, however, it is easier to consider first-passage and last-passage models separately. For completeness, we also state here the first-passage versions of the results above. Define the quantities $\{S(\mathbf{z}), \mathbf{z} \in \mathbb{Z}_+^d\}$ and $\{h(\mathbf{x}), \mathbf{x} \in \mathbb{R}_+^d\}$, analogous to the last-passage quantities $\{T(\mathbf{z})\}$ and $\{g(\mathbf{x})\}$ defined at (2.1) and (2.3), by

(2.9) $$S(\mathbf{z}) = \min_{\pi \in \Pi[\mathbf{z}]} \sum_{\mathbf{v} \in \pi} X(\mathbf{v}),$$

(2.10) $$h(\mathbf{x}) = \inf_{n \in \mathbb{N}} \frac{1}{n} \mathbb{E}\, S(\lfloor n\mathbf{x} \rfloor).$$

COROLLARY 2.5. (i) *If* $\mathbb{E}|X| < \infty$ *and* $\int_{-\infty}^{0} (F(s))^{1/d} < \infty$, *then* $|h(\mathbf{x})| < \infty$ *for all* $\mathbf{x} \in \mathbb{R}_+^d$, *and* $h$ *is continuous on all of* $\mathbb{R}_+^d$ *(including at the boundaries).*

(ii) *Let* $d = 2$. *If* $F$ *has mean* $\mu$ *and variance* $\sigma^2$ *and satisfies* (2.6) *and* (2.7), *then as* $\alpha \downarrow 0$,

$$h(1, \alpha) = \mu - 2\sigma\sqrt{\alpha} + o(\sqrt{\alpha}).$$

Definitions, results and discussions for the growth models are given in Section 5.

## 3. Continuity at the boundary.

3.1. *Case of bounded weights.* We first prove the continuity result of Theorem 2.3 for the case where the weights are bounded. We will need the following concentration inequality, which follows easily from Theorem 8.1.1 of [27]; see, for example, Lemma 5.1 of [20] for the argument.

LEMMA 3.1. *Let* $Y_i, i \in I$, *be a finite collection of independent random variables, such that*

$$\mathbb{P}(|Y_i| \leq L) = 1$$

*for all* $i \in I$. *Let* $\mathcal{C}$ *be a set of subsets of* $I$ *such that*

$$\max_{C \in \mathcal{C}} |C| \leq R,$$

*and let*

$$Z = \max_{C \in \mathcal{C}} \sum_{i \in C} Y_i.$$

*Then for any* $u > 0$,

$$\mathbb{P}(|Z - \mathbb{E}Z| \geq u) \leq \exp\left(-\frac{u^2}{64RL^2} + 64\right).$$

We apply the concentration inequality in the following lemma, which is the central part of the proof of the continuity of $g$.

LEMMA 3.2. *Suppose* $\mathbb{P}(|X| \leq L) = 1$ *for some finite* $L$. *Let* $R > 0$ *and* $\varepsilon > 0$. *There exists* $\delta > 0$ *such that if* $\mathbf{x} \in \mathbb{R}_+^d$ *with* $\|\mathbf{x}\| \leq R$ *and* $x_j = 0$ *(where* $1 \leq j \leq d$), *then*

$$|g(\mathbf{x} + h\mathbf{e}_j) - g(\mathbf{x})| < \varepsilon$$

*for all* $0 \leq h \leq \delta$.

PROOF. Without loss of generality, let $j = 1$. Rephrased, the statement is that for $\mathbf{x} = (x_2, x_3, \ldots, x_d) \in \mathbb{R}_+^{d-1}$,

$$g((h, \mathbf{x})) \to g((0, \mathbf{x})) \qquad \text{as } h \downarrow 0,$$

uniformly in $\{\mathbf{x} : |\mathbf{x}| \le R\}$.

So let $h > 0$ and $n \in \mathbb{N}$. Any path from $\mathbf{0}$ to the point $(\lfloor nh \rfloor, \lfloor n\mathbf{x} \rfloor)$ contains exactly $\lfloor nh \rfloor$ steps which increase the first coordinate, so can be decomposed into a disjoint union of paths from $(r, \mathbf{m}_r)$ to $(r, \mathbf{m}_{r+1})$, $r = 0, 1, 2, \ldots, \lfloor nh \rfloor$, where $\mathbf{m}_r \in \mathbb{Z}_+^{d-1}$ for each $r$, and

$$\mathbf{0} = \mathbf{m}_0 \le \mathbf{m}_1 \le \cdots \le \mathbf{m}_{\lfloor nh \rfloor + 1} = \lfloor n\mathbf{x} \rfloor. \tag{3.1}$$

We have

$$T((\lfloor nh \rfloor, \lfloor n\mathbf{x} \rfloor)) = \max_{\mathbf{m}_0, \mathbf{m}_1, \ldots, \mathbf{m}_{\lfloor nh \rfloor + 1}} \left[ \sum_{r=0}^{\lfloor nh \rfloor} T((r, \mathbf{m}_r), (r, \mathbf{m}_{r+1})) + \sum_{r=0}^{\lfloor nh \rfloor - 1} X(r, \mathbf{m}_{r+1}) \right]. \tag{3.2}$$

(The second term on the right-hand side appears because of the convention that a path from $\mathbf{z}_1$ to $\mathbf{z}_2$ does not include the "final point" $\mathbf{z}_2$.) Here and below the $\mathbf{m}_r$ range over values satisfying (3.1). The number of such choices for the $\mathbf{m}_r$ is

$$\prod_{i=2}^{d} \binom{\lfloor nx_i \rfloor + \lfloor nh \rfloor}{\lfloor nh \rfloor},$$

which, by Stirling's formula, is $\exp[n\phi(h, \mathbf{x}) + o(n)]$, where

$$\phi(h, \mathbf{x}) = \sum_{\substack{2 \le i \le d \\ x_i > 0}} \left( h \log \frac{x_i + h}{h} + x_i \log \frac{x_i + h}{x_i} \right).$$

We now consider the expectation of the quantity within the maximum on the right-hand side of (3.2). For fixed $\{\mathbf{m}_r\}$, we have

$$\mathbb{E}\left[ \sum_{r=0}^{\lfloor nh \rfloor} T((r, \mathbf{m}_r), (r, \mathbf{m}_{r+1})) + \sum_{r=0}^{\lfloor nh \rfloor - 1} X(r, \mathbf{m}_{r+1}) \right]$$

$$= \mathbb{E} \sum_{r=0}^{\lfloor nh \rfloor} T((0, \mathbf{m}_r), (0, \mathbf{m}_{r+1})) + \lfloor nh \rfloor \mathbb{E} X$$

$$\le \mathbb{E} T((0, \mathbf{m}_0), (0, \mathbf{m}_{\lfloor nh \rfloor + 1})) + nhL \qquad \text{(by superadditivity)} \tag{3.3}$$

$$= \mathbb{E} T((0, \lfloor n\mathbf{x} \rfloor)) + nhL$$

$$\le n[g((0, \mathbf{x})) + hL]$$

(by definition of $g$ and superadditivity again).

Still keeping $\{\mathbf{m}_r\}$ fixed, note that the quantity inside the expectation on the left-hand side of (3.3) may be written as the maximum of the sum of various sets of weights $X$; each such set has size $\|(\lfloor nh \rfloor, \lfloor n\mathbf{x} \rfloor)\| \leq n\|(h, \mathbf{x})\|$ and, by assumption, none of the weights has absolute value greater than $L$. So we can apply the concentration inequality in Lemma 3.1 to give

$$\mathbb{P}\left[\sum_{r=0}^{\lfloor nh \rfloor} T((r, \mathbf{m}_r), (r, \mathbf{m}_{r+1})) + \sum_{r=0}^{\lfloor nh \rfloor - 1} X(r, \mathbf{m}_{r+1}) \geq n(g((0, \mathbf{x})) + hL + \varepsilon)\right]$$

$$\leq \mathbb{P}\left[\sum_{r=0}^{\lfloor nh \rfloor} T((r, \mathbf{m}_r), (r, \mathbf{m}_{r+1})) + \sum_{r=0}^{\lfloor nh \rfloor - 1} X(r, \mathbf{m}_{r+1})\right.$$

$$\left. \text{deviates from its expectation by more than } n\varepsilon \right]$$

$$\leq \exp\left(-\frac{(n\varepsilon)^2}{64n\|(h, \mathbf{x})\|L^2} + 64\right)$$

$$= \exp\left(-\frac{n\varepsilon^2}{64\|(h, \mathbf{x})\|L^2} + 64\right).$$

Thus, taking the sum over all possible $\{\mathbf{m}_r\}$,

$$\mathbb{P}[T((\lfloor nh \rfloor, \lfloor n\mathbf{x} \rfloor)) \geq n(g((0, \mathbf{x})) + hL + \varepsilon)]$$

$$\leq \exp[n\phi(h, \mathbf{x}) + o(n)] \exp\left(-\frac{n\varepsilon^2}{64\|(h, \mathbf{x})\|L^2} + 64\right).$$

This sums over $n \in \mathbb{N}$ to a finite amount whenever

$$\varepsilon > 8L\sqrt{\|(h, \mathbf{x})\|\phi(h, \mathbf{x})}.$$

Since

$$g((h, \mathbf{x})) = \lim_{n \to \infty} \frac{1}{n} T((\lfloor nh \rfloor, \lfloor n\mathbf{x} \rfloor)) \qquad \text{a.s.,}$$

Borel–Cantelli then gives

$$g((h, \mathbf{x})) - g((0, \mathbf{x})) \leq hL + 8L\sqrt{\|(h, \mathbf{x})\|\phi(h, \mathbf{x})}.$$

The right-hand side tends to 0 as $h \downarrow 0$, uniformly in $\|\mathbf{x}\| \leq R$, as required.

In the other direction, the superadditivity property in Proposition 2.1(iv) implies

$$g((h, \mathbf{x})) - g((0, \mathbf{x})) \geq g((h, \mathbf{0}))$$
$$= h\mathbb{E} X$$
$$\geq -hL,$$

which again tends to 0 as $h \downarrow 0$, uniformly over all $\mathbf{x}$. □

LEMMA 3.3. *Suppose $\mathbb{P}(|X| \leq L) = 1$ for some finite $L$. Then $g$ is continuous on all of $\mathbb{R}_+^d$.*

PROOF. Let $\varepsilon > 0$ and $\mathbf{y} \in \mathbb{R}_+^d$. Suppose that $\mathbf{y}$ has exactly $k$ nonzero coordinates. Without loss of generality, assume that in fact $y_i > 0$ for $1 \leq i \leq k$ and $y_i = 0$ for $k+1 \leq i \leq d$.

Define the function $g_k$ on $\mathbb{R}_+^k$ by $g_k(\mathbf{u}) = g(u_1, \ldots, u_k, 0, \ldots, 0)$ (appending $d - k$ zeros to the end of $\mathbf{u}$).

Since $g$ is concave on $\mathbb{R}_+^d$ [by Proposition 2.1(iv)], $g_k$ is concave on $\mathbb{R}_+^k$, and so is continuous on the interior of $\mathbb{R}_+^k$.

Hence we can choose $\delta' > 0$ small enough that if $|h_i| < \delta'$ for $1 \leq i \leq k$, then
$$|g_k(x_1 + h_1, \ldots, x_k + h_k) - g_k(x_1, \ldots, x_k)| < \varepsilon,$$
and so

$$\left| g\left(\mathbf{x} + \sum_{i=1}^{k} h_i \mathbf{e}_i\right) - g(\mathbf{x}) \right| < \varepsilon. \tag{3.4}$$

Choose any $R > \|\mathbf{y}\| + k\delta'$. We now fix $\delta > 0$ small enough that the conclusion of Lemma 3.2 applies (for our chosen $L$, $R$ and $\varepsilon$), and also small enough that $\|\mathbf{y}\| + k\delta' + (d-k)\delta < R$.

Take any $\mathbf{h} \in \mathbb{R}^d$ with $\|\mathbf{h}\| \leq \min(\delta', \delta)$ and with $(\mathbf{y} + \mathbf{h}) \in \mathbb{R}_+^d$. Then certainly $|h_i| \leq \delta'$ for $1 \leq i \leq k$, and also $0 \leq h_i \leq \delta$ for $k+1 \leq i \leq d$.

We are about to apply Lemma 3.2 $(d-k)$ times, once for each of the coordinates of $\mathbf{y}$ which is 0. Specifically, for $k+1 \leq j \leq d$, set $\mathbf{x}^{(j)} = \mathbf{y} + \sum_{i=1}^{j-1} h_i \mathbf{e}_i$. Then for $k+1 \leq j \leq d$, we have $\|\mathbf{x}^{(j)}\| < R$ (by choice of $\delta$) and $x_j^{(j)} = 0$ (since $y_j = 0$), so all the required conditions of Lemma 3.2 apply.

Using also (3.4), we obtain

$$|g(\mathbf{y} + \mathbf{h}) - g(\mathbf{y})| \leq \left| g\left(\mathbf{y} + \sum_{i=1}^{k} h_i \mathbf{e}_i\right) - g(\mathbf{y}) \right|$$

$$+ \sum_{j=k+1}^{d} \left| g\left(\mathbf{y} + \sum_{i=1}^{j} h_i \mathbf{e}_i\right) - g\left(\mathbf{y} + \sum_{i=1}^{j-1} h_i \mathbf{e}_i\right) \right|$$

$$= \left| g\left(\mathbf{y} + \sum_{i=1}^{k} h_i \mathbf{e}_i\right) - g(\mathbf{y}) \right|$$

$$+ \sum_{j=k+1}^{d} |g(\mathbf{x}^{(j)} + h_j \mathbf{e}_j) - g(\mathbf{x}^{(j)})|$$

$$< \varepsilon + (d-k)\varepsilon$$

$$\leq (d+1)\varepsilon.$$

Since $\varepsilon$ was arbitrary, we have that $g$ is continuous at $\mathbf{y}$, as desired. $\square$

3.2. *Extension to unbounded weight distribution.* Define a lattice animal of size $n$ to be a connected subset of $\mathbb{Z}^d$ of size $n$ which includes the origin. Extend the i.i.d. array $\{X(\mathbf{z})\}$ to all of $\mathbb{Z}^d$, so that $\{X(\mathbf{z}), \mathbf{z} \in \mathbb{Z}^d\}$ is an i.i.d. array with common distribution $F$. Let $A(n)$ be the set of lattice animals of size $n$, and define

$$N(n) = \max_{\xi \in A(n)} \sum_{\mathbf{z} \in \xi} X(\mathbf{z}),$$

the maximum weight of a lattice animal of size $n$.

The results in the following proposition are taken from Theorems 1.1 and 2.3 of [20]. (In fact the model in [20] covers only the case where the weights $X(\mathbf{z})$ are nonnegative; but replacing weights whose value is 0 by negative weights can only reduce the left-hand side of (3.5) or (3.6) and leaves the right-hand side unchanged. Alternatively, see [6] for a detailed treatment of the lattice animals model where the weights can take negative values.)

PROPOSITION 3.4. *There exists $c = c(d) < \infty$ such that, for all $F$ satisfying* (2.5):

(i) *for all $n > 1$,*

(3.5) $$\mathbb{E} N(n) \leq cn \int_0^\infty (1 - F(s))^{1/d} \, ds;$$

(ii) *with probability 1,*

(3.6) $$\limsup_{n \to \infty} \frac{N(n)}{n} \leq c \int_0^\infty (1 - F(s))^{1/d} \, ds.$$

Now we can easily deduce the following lemma for the directed percolation model. Part (iii) implies Proposition 2.2.

LEMMA 3.5. *There exists $c = c(d) < \infty$ such that, for all $F$ satisfying* (2.5):

(i) *for all $\mathbf{z} \in \mathbb{Z}_+^d$,*

(3.7) $$\mathbb{E} T(\mathbf{z}) \leq c\|\mathbf{z}\| \int_0^\infty (1 - F(s))^{1/d} \, ds;$$

(ii) *with probability 1,*

(3.8) $$\limsup_{n \to \infty} \frac{1}{n} \max_{\mathbf{z} : \|\mathbf{z}\| \leq n} T(\mathbf{z}) \leq c \int_0^\infty (1 - F(s))^{1/d} \, ds;$$

(iii) *for all* $\mathbf{x} \in \mathbb{R}_+^d$,

(3.9) $$\|\mathbf{x}\| \mathbb{E} X \leq g(\mathbf{x}) \leq c \|\mathbf{x}\| \int_0^\infty (1 - F(s))^{1/d} \, ds.$$

PROOF. First note that if $\mathbf{z} \in \mathbb{Z}_+^d$ and $\|\mathbf{z}\| = n$, then any path $\pi \in \Pi[\mathbf{z}]$ is a lattice animal of size $n$; thus $T(\mathbf{z}) \leq N(n)$, and parts (i) and (ii) follow immediately from Proposition 3.4.

Putting $\mathbf{z} = \lfloor n\mathbf{x} \rfloor$ in (i), dividing by $n$ and then letting $n \to \infty$ gives the upper bound in (iii).

For the lower bound in (iii), let $\mathbf{z} \in \mathbb{Z}_+^d$, and let $\tilde{\pi}$ be any path in $\Pi[\mathbf{z}]$; then $|\tilde{\pi}| = \|\mathbf{z}\|$, and we have

$$\mathbb{E} T(\mathbf{z}) = \mathbb{E} \max_{\pi \in \Pi[\mathbf{z}]} \sum_{\mathbf{v} \in \pi} X(\mathbf{v})$$
$$\geq \mathbb{E} \sum_{\mathbf{v} \in \tilde{\pi}} X(\mathbf{v})$$
$$= \|\mathbf{z}\| \mathbb{E} X.$$

Again let $\mathbf{z} = \lfloor n\mathbf{x} \rfloor$ and let $n \to \infty$ to obtain the lower bound in (iii). □

We now introduce truncated versions of the weights $\{X(\mathbf{z})\}$. For $L > 0$ and $\mathbf{z} \in \mathbb{Z}_+^d$, let $X^{(L)}(\mathbf{z}) = \max\{\min\{X(\mathbf{z}), L\}, -L\}$ [so that $|X^{(L)}(\mathbf{z})| = \min(|X(\mathbf{z})|, L)$].

Then let $\{T^{(L)}(\mathbf{z}), \mathbf{z} \in \mathbb{Z}_+^d\}$ and $\{g^{(L)}(\mathbf{x}), \mathbf{x} \in \mathbb{R}_+^d\}$ be defined just as $\{T(\mathbf{z})\}$ and $\{g(\mathbf{x})\}$, but with the quantities $\{X(\mathbf{z})\}$ replaced by the truncated versions $\{X^{(L)}(\mathbf{z})\}$.

LEMMA 3.6. *Suppose that* (2.4) *and* (2.5) *hold. Then for any* $\mathbf{x} \in \mathbb{R}_+^d$,

(3.10) $$g^{(L)}(\mathbf{x}) - \|\mathbf{x}\| \int_{-\infty}^{-L} F(s) \, ds$$
$$\leq g(\mathbf{x}) \leq g^{(L)}(\mathbf{x}) + c\|\mathbf{x}\| \int_L^\infty (1 - F(s))^{1/d} \, ds,$$

*where* $c$ *is as in Lemma* 3.5. *Thus, for any* $R > 0$,

(3.11) $$\sup_{\mathbf{x} \in \mathbb{R}_+^d : \|\mathbf{x}\| \leq R} |g(\mathbf{x}) - g^{(L)}(\mathbf{x})| \to 0 \qquad \text{as } L \to \infty.$$

PROOF. Note that

$$-[L - X(\mathbf{z})]_+ \leq X(\mathbf{z}) - X^{(L)}(\mathbf{z}) \leq [X(\mathbf{z}) - L]_+.$$

We consider first the positive tail. Let $\mathbf{x} \in \mathbb{R}_+^d$. Then

$$g(\mathbf{x}) - g^{(L)}(\mathbf{x}) = \lim_{n\to\infty} \frac{1}{n}\mathbb{E}\,T(\lfloor n\mathbf{x}\rfloor) - \lim_{n\to\infty} \frac{1}{n}\mathbb{E}\,T^{(L)}(\lfloor n\mathbf{x}\rfloor)$$

$$= \lim_{n\to\infty} \frac{1}{n}\mathbb{E}\left[\sup_{\pi\in\Pi[\mathbf{0},\lfloor n\mathbf{x}\rfloor)}\sum_{\mathbf{v}\in\pi} X(\mathbf{v}) - \sup_{\pi\in\Pi[\mathbf{0},\lfloor n\mathbf{x}\rfloor)}\sum_{\mathbf{v}\in\pi} X^{(L)}(\mathbf{v})\right]$$

(3.12)
$$\leq \lim_{n\to\infty} \frac{1}{n}\mathbb{E}\sup_{\pi\in\Pi[\mathbf{0},\lfloor n\mathbf{x}\rfloor)}\left[\sum_{\mathbf{v}\in\pi} X(\mathbf{v}) - \sum_{\mathbf{v}\in\pi} X^{(L)}(\mathbf{v})\right]$$

$$= \lim_{n\to\infty} \frac{1}{n}\mathbb{E}\sup_{\pi\in\Pi[\mathbf{0},\lfloor n\mathbf{x}\rfloor)}\left[\sum_{\mathbf{v}\in\pi}[X(\mathbf{v}) - L]_+\right]$$

$$\leq c\|\mathbf{x}\|\int_L^\infty (1-F(s))^{1/d}\,ds;$$

the last inequality follows from Lemma 3.5, since the variables $\{[X(\mathbf{v}) - L]_+, \mathbf{v} \in \mathbb{Z}_+^d\}$ are i.i.d. with distribution function $F^{(>L)}$, where $F^{(>L)}(s) = 0$ for $s \leq L$ and $F^{(>L)}(s) = F(s)$ for $s > L$. This gives the lower bound in (3.10).

For the negative tail, let $\mathbf{z} \in \mathbb{Z}_+^d$, and let $\pi^* \in \Pi[\mathbf{0}, \mathbf{z}]$ be the maximizing path for $T^{(L)}(\mathbf{z})$. (If there are several maximizing paths, choose, say, the one that is first in the lexicographic order.)

Now for any $\mathbf{v} \in \mathbb{Z}_+^d$, $\mathbb{P}(\mathbf{v} \in \pi^* | X(\mathbf{v}) \leq s)$ is a nondecreasing function of $s$. [This follows from a simple coupling, since $\{X(\mathbf{v}'), \mathbf{v}' \in \mathbb{Z}_+^d\}$ are independent and, for a fixed realization of the other weights $\{X(\mathbf{v}'), \mathbf{v}' \in \mathbb{Z}_+^d, \mathbf{v}' \neq \mathbf{v}\}$, the function $I\{\mathbf{v} \in \pi^*\}$ is a nondecreasing function of $X(\mathbf{v})$.]

Hence in particular $\mathbb{P}(\mathbf{v} \in \pi^* | X(\mathbf{v}) \leq -L) \leq \mathbb{P}(\mathbf{v} \in \pi^*)$, and so by simple manipulation of conditional probabilities,

$$\mathbb{P}(X(\mathbf{v}) \leq -L | \mathbf{v} \in \pi^*) \leq \mathbb{P}(X(\mathbf{v}) \leq -L).$$

Furthermore, the event $\{\mathbf{v} \in \pi^*\}$ depends only on $\{\max\{X(\mathbf{v}'), -L\}, \mathbf{v}' \in \mathbb{Z}_+^d\}$; thus, conditional on $\{X(\mathbf{v}) \leq -L\}$, $X(\mathbf{v})$ is independent of $\{\mathbf{v} \in \pi^*\}$. Hence

$$\mathbb{E}\left([-L - X(\mathbf{v})]_+ | \mathbf{v} \in \pi^*\right)$$
$$= \mathbb{E}\left([-L - X(\mathbf{v})]_+ | X(\mathbf{v}) \leq -L\right)\mathbb{P}(X(\mathbf{v}) \leq -L | \mathbf{v} \in \pi^*)$$
$$\leq \mathbb{E}\left([-L - X(\mathbf{v})]_+ | X(\mathbf{v}) \leq -L\right)\mathbb{P}(X(\mathbf{v}) \leq -L)$$
$$= \mathbb{E}\left([-L - X(\mathbf{v})]_+\right)$$
$$= \int_{-\infty}^{-L} F(s)\,ds.$$

Now

$$\mathbb{E}\, T(\mathbf{z}) = \mathbb{E} \max_{\pi \in \Pi[\mathbf{z}]} \sum_{\mathbf{v} \in \pi} X(\mathbf{v})$$

$$\geq \mathbb{E} \max_{\pi \in \Pi[\mathbf{z}]} \sum_{\mathbf{v} \in \pi} (X^{(L)}(\mathbf{v}) - [-L - X(\mathbf{v})]_+)$$

$$\geq \mathbb{E} \sum_{\mathbf{v} \in \pi^*} X^{(L)}(\mathbf{v}) - \mathbb{E} \sum_{\mathbf{v} \in \pi^*} [-L - X(\mathbf{v})]_+$$

$$= \mathbb{E}\, T^{(L)}(\mathbf{z}) - \sum_{\mathbf{v} \in \mathbb{Z}_+^d} \mathbb{P}(\mathbf{v} \in \pi^*) \mathbb{E}\left([-L - X(\mathbf{v})]_+ | \mathbf{v} \in \pi^*\right)$$

$$\geq \mathbb{E}\, T^{(L)}(\mathbf{z}) - \int_{-\infty}^{-L} F(s)\, ds \sum_{\mathbf{v} \in \mathbb{Z}_+^d} \mathbb{P}(\mathbf{v} \in \pi^*)$$

$$= \mathbb{E}\, T^{(L)}(\mathbf{z}) - \int_{-\infty}^{-L} F(s)\, ds\, \mathbb{E}\, |\pi^*|$$

$$= \mathbb{E}\, T^{(L)}(\mathbf{z}) - \|\mathbf{z}\| \int_{-\infty}^{-L} F(s)\, ds,$$

since $\pi^* \in \Pi[\mathbf{z}]$ and all paths in $\Pi[\mathbf{z}]$ have length $\|\mathbf{z}\|$.

Letting $\mathbf{z} = \lfloor n\mathbf{x} \rfloor$ and taking $n \to \infty$ then gives the first inequality in (3.10), as required. The convergence in (3.11) then follows since, under (2.4) and (2.5), both the integrals in (3.10) tend to 0 as $L \to \infty$. $\square$

It is now immediate to extend the continuity property proved in Lemma 3.3 to the case of unbounded weights, and so complete the proof of Theorem 2.3.

Lemma 3.3 shows that the functions $g^{(L)}$ are continuous for each $L$, and (3.11) shows that, under the conditions (2.4) and (2.5), $g^{(L)} \to g$ as $L \to \infty$, uniformly on any compact subset of $\mathbb{R}_+^d$. Hence $g$ itself is continuous, as required.

**4. Asymptotics at the boundary for $d = 2$.** In this section we prove Theorem 2.4.

We first obtain an estimate on the growth rate in the case where $F$ is a Bernoulli distribution. This is done using a comparison with an alternative percolation model in which the Bernoulli distribution is an exactly solvable case. We write $\mathrm{Ber}(p)$ for the Bernoulli distribution with parameter $p$, with $\mathbb{P}(X = 1) = 1 - \mathbb{P}(X = 0) = p$.

LEMMA 4.1. *For all $\alpha > 0$, $p \in [0, 1]$,*

$$g_{\mathrm{Ber}(p)}(1, \alpha) \leq (1 + \alpha)p + 2\sqrt{\alpha}\sqrt{1 + \alpha}\sqrt{p(1-p)}.$$

PROOF. Recall that
$$T(m,n) = \max_{\pi \in \Pi[m,n]} \sum_{\mathbf{z} \in \pi} X(\mathbf{z}),$$
where $\Pi[m,n]$ is the set of paths of the form
$$\mathbf{z}_0, \mathbf{z}_1, \ldots, \mathbf{z}_{m+n-1}$$
such that $\mathbf{z}_0 = \mathbf{0}$, such that, for all $1 \leq i \leq m+n-1$, $\mathbf{z}_i - \mathbf{z}_{i-1} = \mathbf{e}_j$ for some $j \in \{1,\ldots,d\}$, and such that also $(m,n) - \mathbf{z}_{m+n-1} = \mathbf{e}_j$ for some $j \in \{1,\ldots,d\}$.

Define an alternative set of increasing paths $\widetilde{\Pi}[m,n]$ to be those paths of the form
$$(0, y_0), (1, y_1), \ldots, (m-1, y_{m-1}),$$
where $0 \leq y_0 \leq y_1 \leq \cdots \leq y_{m-1} \leq n$, and define
$$\widetilde{T}(m,n) = \max_{\tilde{\pi} \in \widetilde{\Pi}[m,n]} \sum_{\mathbf{z} \in \tilde{\pi}} X(\mathbf{z}).$$

Define the function $\psi: \mathbb{Z}_+^2 \mapsto \mathbb{Z}_+^2$ by $\psi(x,y) = (x+y, y)$.

Since $\psi(\mathbf{z}) \neq \psi(\mathbf{z}')$ whenever $\mathbf{z} \neq \mathbf{z}'$, and since $\{X(\mathbf{z}), \mathbf{z} \in \mathbb{Z}_+^2\}$ are i.i.d., we have

(4.1) $\quad \mathbb{E}\, T(m,n) = \mathbb{E} \max_{\pi \in \Pi[m,n]} \sum_{\mathbf{z} \in \pi} X(\mathbf{z}) = \mathbb{E} \max_{\pi \in \Pi[m,n]} \sum_{\mathbf{z} \in \pi} X(\psi(\mathbf{z})).$

For a path $\pi = \mathbf{z}_0, \ldots, \mathbf{z}_{r-1}$, write $\psi(\pi)$ for the path $\psi(\mathbf{z}_0), \ldots, \psi(\mathbf{z}_{r-1})$. From the definitions of the path sets above, one can obtain that if $\pi \in \Pi[m,n]$, then $\psi(\pi) \in \widetilde{\Pi}[m+n, n]$. (We do not require that all paths in $\widetilde{\Pi}[m+n,n]$ can be written as $\psi(\pi)$ in this way.) Put another way: let $\psi(\Pi[m,n])$ be the set $\{\tilde{\pi}: \tilde{\pi} = \psi(\pi)$ for some $\pi \in \Pi[m,n]\}$; then $\psi(\Pi[m,n]) \subseteq \widetilde{\Pi}[m+n,n]$.

Continuing from (4.1),
$$\mathbb{E}\, T(m,n) = \mathbb{E} \max_{\pi \in \Pi[m,n]} \sum_{\mathbf{v} \in \psi(\pi)} X(\mathbf{v})$$
$$= \mathbb{E} \max_{\tilde{\pi} \in \psi(\Pi[m,n])} \sum_{\mathbf{v} \in \tilde{\pi}} X(\mathbf{v})$$
$$\leq \mathbb{E} \max_{\tilde{\pi} \in \widetilde{\Pi}[m+n,n]} \sum_{\mathbf{v} \in \tilde{\pi}} X(\mathbf{v})$$
$$= \mathbb{E}\, \widetilde{T}(m+n, n).$$

Then

(4.2)
$$g(1, \alpha) = \lim_{n \to \infty} \frac{1}{n} \mathbb{E}\, T(n, \lfloor \alpha n \rfloor)$$
$$\leq \liminf_{n \to \infty} \frac{1}{n} \mathbb{E}\, \widetilde{T}(\lfloor 1+\alpha \rfloor n, \lfloor \alpha n \rfloor).$$

Seppäläinen [26] analyzes directed percolation based on the path sets $\widetilde{\Pi}$, and in particular obtains that, for the case of Bernoulli weights,

$$\lim_{n\to\infty}\frac{1}{n}\mathbb{E}_{\mathrm{Ber}(p)}\widetilde{T}(\lfloor\alpha_1 n\rfloor,\lfloor\alpha_2 n\rfloor)$$
$$=\begin{cases} p(\alpha_1-\alpha_2)+2\sqrt{\alpha_1\alpha_2}\sqrt{p(1-p)}, & \text{if } p\leq \dfrac{\alpha_1}{\alpha_1+\alpha_2},\\ \alpha_1, & \text{if } p\geq \dfrac{\alpha_1}{\alpha_1+\alpha_2}.\end{cases}$$

A calculation then shows that for all $p$,

$$\lim_{n\to\infty}\frac{1}{n}\mathbb{E}_{\mathrm{Ber}(p)}\widetilde{T}(\lfloor\alpha_1 n\rfloor,\lfloor\alpha_2 n\rfloor)\leq \alpha_1 p+2\sqrt{\alpha_1\alpha_2}\sqrt{p(1-p)}.$$

Substituting into (4.2) with $\alpha_1=1+\alpha$, $\alpha_2=\alpha$ gives the required result. $\square$

LEMMA 4.2. *Let $F_1$ and $F_2$ be distributions with means $\mu_1$ and $\mu_2$ and satisfying (2.6) and (2.7). Then for all $\alpha>0$,*

$$|g_{F_1}(1,\alpha)-g_{F_2}(1,\alpha)-(1+\alpha)(\mu_1-\mu_2)|$$
$$\leq 2\sqrt{\alpha(1+\alpha)}\int_{-\infty}^{\infty}|F_1(s)-F_2(s)|^{1/2}\,ds.$$

PROOF. Let $\{U(\mathbf{z}),\mathbf{z}\in\mathbb{Z}_+^2\}$ be i.i.d. uniform on $[0,1]$, and for $i=1,2$, let $X_i(\mathbf{z})=F_i^{-1}(U(\mathbf{z}))$, where $F^{-1}(u)=\sup\{x:F(x)\leq u\}$.

Then for $i=1,2$, $\{X_i(\mathbf{z}),\mathbf{z}\in\mathbb{Z}_+^2\}$ are i.i.d. with distribution $F_i$, and for any $x$,

$$\mathbb{P}(X_1(\mathbf{z})\geq x, X_2(\mathbf{z})<x)=[F_2(x)-F_1(x)]_+$$

and

$$\mathbb{P}(X_2(\mathbf{z})\geq x, X_1(\mathbf{z})<x)=[F_1(x)-F_2(x)]_+.$$

Now

$$g_{F_1}(1,\alpha)-g_{F_2}(1,\alpha)$$
$$=\lim_{n\to\infty}\frac{1}{n}\mathbb{E}\max_{\pi\in\Pi(n,\lfloor\alpha n\rfloor)}\sum_{\mathbf{z}\in\pi}X_1(\mathbf{z})-\lim_{n\to\infty}\mathbb{E}\max_{\pi\in\Pi(n,\lfloor\alpha n\rfloor)}\sum_{\mathbf{z}\in\pi}X_2(\mathbf{z})$$
$$\leq \lim_{n\to\infty}\frac{1}{n}\mathbb{E}\max_{\pi\in\Pi(n,\lfloor\alpha n\rfloor)}\sum_{\mathbf{z}\in\pi}(X_1(\mathbf{z})-X_2(\mathbf{z}))$$
$$=\lim_{n\to\infty}\frac{1}{n}\mathbb{E}\max_{\pi\in\Pi(n,\lfloor\alpha n\rfloor)}\int_{-\infty}^{\infty}\sum_{\mathbf{z}\in\pi}[I(X_1(\mathbf{z})\geq x, X_2(\mathbf{z})<x)$$
$$-I(X_1(\mathbf{z})<x, X_2(\mathbf{z})\geq x)]\,dx$$

$$\leq \lim_{n\to\infty} \frac{1}{n} \mathbb{E} \int_{-\infty}^{\infty} \max_{\pi\in\Pi(n,\lfloor\alpha n\rfloor)} \sum_{\mathbf{z}\in\pi} [I(X_1(\mathbf{z}) \geq x, X_2(\mathbf{z}) < x)$$
$$- I(X_1(\mathbf{z}) < x, X_2(\mathbf{z}) \geq x)] \, dx$$
$$= \int_{-\infty}^{\infty} \lim_{n\to\infty} \frac{1}{n} \mathbb{E} \max_{\pi\in\Pi(n,\lfloor\alpha n\rfloor)} \sum_{\mathbf{z}\in\pi} [I(X_1(\mathbf{z}) \geq x, X_2(\mathbf{z}) < x)$$
$$- I(X_1(\mathbf{z}) < x, X_2(\mathbf{z}) \geq x)] \, dx$$

(by Fubini's theorem and bounded convergence)

$$\leq \int_{-\infty}^{\infty} \lim_{n\to\infty} \frac{1}{n} \left\{ \mathbb{E} \max_{\pi\in\Pi(n,\lfloor\alpha n\rfloor)} \sum_{\mathbf{z}\in\pi} [I(X_1(\mathbf{z}) \geq x, X_2(\mathbf{z}) < x)] \right.$$
$$\left. + \mathbb{E} \max_{\pi\in\Pi(n,\lfloor\alpha n\rfloor t)} \left( \sum_{\mathbf{z}\in\pi} [1 - I(X_1(\mathbf{z}) < x, X_2(\mathbf{z}) \geq x)] \right. \right.$$
$$\left. \left. - \sum_{\mathbf{z}\in\pi} 1 \right) \right\} dx$$
$$= \int_{-\infty}^{\infty} \{ g_{\mathrm{Ber}([F_2(x)-F_1(x)]_+)}(1,\alpha) + g_{\mathrm{Ber}(1-[F_1(x)-F_2(x)]_+)}(1,\alpha) - (1+\alpha) \, dx \}$$
$$\leq \int_{-\infty}^{\infty} \{ (1+\alpha)([F_2(x)-F_1(x)]_+ + 1 - [F_1(x)-F_2(x)]_+)$$
$$+ 2\sqrt{\alpha}\sqrt{1+\alpha}([F_2(x)-F_1(x)]_+^{1/2} + [F_1(x)-F_2(x)]_+^{1/2})$$
$$- (1+\alpha) \} \, dx$$

(by Lemma 4.1)

$$= (1+\alpha) \int_{-\infty}^{\infty} \{F_2(x) - F_1(x)\} \, dx$$
$$+ 2\sqrt{\alpha}\sqrt{1+\alpha} \int_{-\infty}^{\infty} |F_1(x) - F_2(x)|^{1/2} \, dx$$
$$= (1+\alpha)(\mu_1 - \mu_2) + 2\sqrt{\alpha}\sqrt{1+\alpha} \int_{-\infty}^{\infty} |F_1(x) - F_2(x)|^{1/2} \, dx.$$

Similarly,

$$g_{F_2}(1,\alpha) - g_{F_1}(1,\alpha)$$
$$\leq (1+\alpha)(\mu_2 - \mu_1) + 2\sqrt{\alpha}\sqrt{1+\alpha} \int_{-\infty}^{\infty} |F_1(x) - F_2(x)|^{1/2} \, dx.$$

Together these give the desired result. □

LEMMA 4.3. *Let $F$ satisfy (2.6) and (2.7), and let $\varepsilon > 0$. Then there is a distribution $\widetilde{F}$ with bounded support which has the same mean and variance as $F$, and which satisfies*

$$\int_{-\infty}^{\infty} |F(s) - \widetilde{F}(s)|^{1/2} \, ds < \varepsilon.$$

PROOF. Let $X$ have distribution $F$. For brevity we cover only the case where $\mathbb{P}(X \geq 0) = 1$; the negative tail can be truncated in an analogous way.

Take any $t > 0$. If $\mathbb{P}(X > t) = 0$, then $F$ itself has bounded support and we take $\widetilde{F} = F$. Otherwise, let $m = \mathbb{E}(X|X > t)$ and $w = \mathbb{E}(X^2|X > t)$, and choose $p$, $u$ to satisfy

(4.3) $$(1-p)t + pu = m,$$

(4.4) $$(1-p)t^2 + pu^2 = w.$$

The solution is

$$p = \frac{(m-t)^2}{(m-t)^2 + w - m^2},$$

$$u = t + \frac{m-t}{p};$$

note that $0 < p \leq 1$ and $u > t$, since $m > t$ and $w \geq m^2$.

Then let

(4.5) $$\widetilde{F}(x) = \begin{cases} F(x), & \text{if } 0 \leq x < t, \\ 1 - p[1 - F(t)], & \text{if } t \leq x < u, \\ 1, & \text{if } x \geq u. \end{cases}$$

Now $\widetilde{F}$ has bounded support since $\widetilde{F}(u) = 1$. For the mean and variance, let $\widetilde{X}$ have distribution $\widetilde{F}$, and use (4.3) and (4.5) to give

$$\mathbb{E}\widetilde{X} = \mathbb{E}(X; X \leq t) + (1-p)[1 - F(t)]t + p[1 - F(t)]u$$
$$= \mathbb{E}(X; X \leq t) + \mathbb{P}(X > t)[(1-p)t + pu]$$
$$= \mathbb{E}X;$$

similarly use (4.4) to give $\mathbb{E}\widetilde{X}^2 = \mathbb{E}X^2$.

For the final part we have

$$\int_0^\infty |F(s) - \widetilde{F}(s)|^{1/2} \, ds \leq \int_t^\infty [1 - F(s)]^{1/2} \, ds + \int_t^\infty [1 - \widetilde{F}(s)]^{1/2} \, ds$$
$$= \int_t^\infty [1 - F(s)]^{1/2} \, ds + \int_t^u [p(1 - F(t))]^{1/2} \, ds$$
$$\leq \int_t^\infty [1 - F(s)]^{1/2} \, ds + [pu^2 \mathbb{P}(X > t)]^{1/2}$$
$$\leq \int_t^\infty [1 - F(s)]^{1/2} \, ds + \mathbb{E}(X^2; X > t)^{1/2}.$$

By assumption, $\int_0^\infty [1 - F(s)]^{1/2}\, ds < \infty$, and this implies that $\mathbb{E}\, X^2 < \infty$ also; hence by choosing $t$ large enough we can make the right-hand side as small as desired. $\square$

LEMMA 4.4. *Let $F$ be a distribution with bounded support, and, for $k \in \mathbb{N}$, let $F^{(k)}$ be the distribution of $X_1 + X_2 + \cdots + X_k$, where $\{X_i\}$ are i.i.d. $\sim F$. Let $r:\mathbb{R}_+ \mapsto \mathbb{N}$ be any function satisfying $r(\alpha) \to \infty$ and $r(\alpha)\sqrt{\alpha} \to 0$ as $\alpha \downarrow 0$. Then*

$$\lim_{\alpha \downarrow 0} \frac{1}{\sqrt{\alpha}}\left| g_F(1,\alpha) - \frac{1}{r(\alpha)} g_{F^{(r)}}(1, \alpha r(\alpha)) \right| = 0.$$

PROOF. For $(x,y) \in \mathbb{Z}_+^2$ and $r \in \mathbb{N}$, let $B^{(r)}(x,y)$ be the set $\{(rx + i, y), i = 0, 1, \ldots, r-1\}$.

The sets $B^{(r)}(\mathbf{z})$ partition $\mathbb{Z}_+^2$; essentially we have grouped the sites of $\mathbb{Z}_+^2$ into "blocks" of length $r$ and height 1. We will compare our original model with one where each of these blocks functions as a single site, whose weight is the sum of the original sites contained in the block.

Given $\pi \in \Pi(nr, m)$, we can find $\tilde{\pi} \in \Pi(n, m)$ such that

(4.6) $$\left| \bigcup_{\mathbf{z} \in \tilde{\pi}} B^{(r)}(\mathbf{z}) \triangle \pi \right| \leq mr,$$

where $\triangle$ denotes the symmetric difference. [For instance, let $\tilde{\pi} = \{\mathbf{z} : \pi \cap B^{(r)}(\mathbf{z}) \neq \varnothing\}$.] Similarly, given $\tilde{\pi} \in \Pi(n, m)$, we can find $\pi \in \Pi(nr, m)$ such that (4.6) is satisfied.

Suppose that $\{X(\mathbf{z}), \mathbf{z} \in \mathbb{Z}_+^2\}$ are i.i.d. with distribution $F$. Define

$$\widetilde{X}^{(r)}(\mathbf{z}) = \sum_{\mathbf{v} \in B^{(r)}(\mathbf{z})} X(\mathbf{z}).$$

Then $\{\widetilde{X}^{(r)}(\mathbf{z}), \mathbf{z} \in \mathbb{Z}_+^2\}$ are i.i.d. with distribution $F^{(r)}$. Let $K$ be such that $\mathbb{P}(|X| > K) = 0$. Then by the properties (in both directions) noted at (4.6) and after,

$$\left| \max_{\pi \in \Pi(nr,m)} \sum_{\mathbf{z} \in \pi} X(\mathbf{z}) - \max_{\tilde{\pi} \in \Pi(n,m)} \sum_{\mathbf{z} \in \tilde{\pi}} \widetilde{X}(\mathbf{z}) \right| \leq mrK,$$

so that

$$\left| \frac{1}{nr} \mathbb{E}_F T(nr, m) - \frac{1}{r}\frac{1}{n} \mathbb{E}_{F^{(r)}} T(n, m) \right| \leq \frac{1}{nr} mrK.$$

Putting $m = \lfloor \alpha nr \rfloor$ and letting $n \to \infty$ gives

$$\left| g_F(1, \alpha) - \frac{1}{r} g_{F^{(r)}}(1, \alpha r) \right| \leq \alpha r K.$$

If $r$ is a function of $\alpha$ such that $r\sqrt{\alpha} \to 0$ as $\alpha \downarrow 0$, then the right-hand side is $o(\sqrt{\alpha})$ as $\alpha \downarrow 0$, as desired. $\square$

LEMMA 4.5. *Let $F$ be a distribution with bounded support, with mean $\mu_F$ and variance $\sigma_F^2$. Let $F^{(k)}$ and $r$ be as in Lemma* 4.4. *Then*

$$\lim_{\alpha \downarrow 0} \frac{1}{\sqrt{\alpha}} \left| \frac{1}{r(\alpha)} g_{F^{(r)}}(1, \alpha r(\alpha)) - \mu_F - \sigma_F \frac{g_\Phi(1, \alpha r(\alpha))}{\sqrt{r(\alpha)}} \right| = 0,$$

*where $\Phi$ is the standard normal distribution.*

PROOF. Theorem 5.16 of [24] gives a bound on the rate of convergence in the central limit theorem, for distributions $F$ which have a finite third moment. Here $F$ has bounded support and so certainly finite third moment; we obtain that there exists $C = C(F)$ such that

$$|\widetilde{F}^{(r)}(x) - \Phi(x)| \leq C r^{-1/2} (1 + |x|)^{-3}$$

for all $r \in \mathbb{N}$, $x \in \mathbb{R}$, where $\widetilde{F}^{(r)}$ is the distribution of $(X_1 + \cdots + X_r - r\mu_F)/(\sigma_F \sqrt{r})$. (Note that $\widetilde{F}^{(r)}$ is simply the distribution $F^{(r)}$ normalized to have mean 0 and variance 1.)

Now combine this estimate with Lemma 4.2; for any $r \in \mathbb{N}$,

$$\frac{1}{\sqrt{\alpha}} \left| \frac{1}{r} g_{F^{(r)}}(1, \alpha r) - \mu_F - \sigma_F \frac{g_\Phi(1, \alpha r)}{\sqrt{r}} \right|$$

$$= \frac{1}{\sqrt{\alpha}} \left| \frac{\sigma_F}{\sqrt{r}} g_{\widetilde{F}^{(r)}}(1, \alpha r) - \sigma_F \frac{g_\Phi(1, \alpha r)}{\sqrt{r}} \right|$$

$$= \frac{\sigma_F}{\sqrt{\alpha r}} |g_{\widetilde{F}^{(r)}}(1, \alpha r) - g_\Phi(1, \alpha r)|$$

$$\leq \frac{\sigma_F}{\sqrt{\alpha r}} 2\sqrt{\alpha r} \sqrt{1 + \alpha r} \int_{-\infty}^{\infty} |C r^{-1/2} (1 + |x|)^{-3}|^{1/2} \, dx$$

$$= C' \sigma_F \frac{\sqrt{1 + \alpha r}}{r^{1/4}},$$

where $C'$ is some constant independent of $r$ and $\alpha$. If $r$ is a function of $\alpha$ such that $r\sqrt{\alpha} \to 0$ and $r \to \infty$ as $\alpha \downarrow 0$, then the right-hand side tends to 0 as $\alpha \downarrow 0$, as required. $\square$

The following lemma is the universality result which we need:

LEMMA 4.6. *Let $F$ be a distribution with mean $\mu_F$ and variance $\sigma_F^2$, and satisfying* (2.6) *and* (2.7), *and let the function $r$ satisfy the conditions of Lemma* 4.4. *Then*

(4.7) $$\lim_{\alpha \downarrow 0} \frac{1}{\sqrt{\alpha}} \left| g_F(1, \alpha) - \mu_F - \sigma_F \frac{g_\Phi(1, \alpha r(\alpha))}{\sqrt{r(\alpha)}} \right| = 0.$$

PROOF. Let $\varepsilon > 0$. Using Lemma 4.3, choose a distribution $\widetilde{F}$ with bounded support, with the same mean and variance as $F$, and with

(4.8) $$\int_{-\infty}^{\infty} |F(x) - \widetilde{F}(x)|^{1/2}\, dx < \varepsilon/2.$$

Then

$$\limsup_{\alpha \downarrow 0} \frac{1}{\sqrt{\alpha}} \left| g_F(1,\alpha) - \mu_F - \sigma_F \frac{g_\Phi(1,\alpha r(\alpha))}{\sqrt{r(\alpha)}} \right|$$

$$\leq \limsup_{\alpha \downarrow 0} \frac{1}{\sqrt{\alpha}} \left| g_{\widetilde{F}}(1,\alpha) - \mu_F - \sigma_F \frac{g_\Phi(1,\alpha r(\alpha))}{\sqrt{r(\alpha)}} \right|$$

$$+ \limsup_{\alpha \downarrow 0} \frac{1}{\sqrt{\alpha}} |g_F(1,\alpha) - g_{\widetilde{F}}(1,\alpha)|.$$

The first term is 0 by combining Lemmas 4.4 and 4.5 and using the fact that $\widetilde{F}$ has the same mean and variance as $F$; the second term is $\leq \varepsilon$ by Lemma 4.2 and (4.8). This works for any $\varepsilon > 0$, so the desired result follows. □

Finally, we compare with an exactly solvable case to yield the asymptotic behavior for all $F$:

PROOF OF THEOREM 2.4. Choose any $r$ that satisfies the conditions of Lemma 4.4, for example, $r(\alpha) = \lfloor \alpha^{-1/4} \rfloor$.

When $F$ is the exponential distribution with mean 1 (and so also variance 1), we have the exact formula $g_F(1,\alpha) = 1 + 2\sqrt{\alpha} + \alpha$ (see, e.g., [25]). Substituting into (4.7) gives

(4.9) $$\lim_{\alpha \downarrow 0} \frac{1}{\sqrt{\alpha}} \left| 2\sqrt{\alpha} - \frac{g_\Phi(1,\alpha r(\alpha))}{\sqrt{r(\alpha)}} \right| = 0.$$

Now take any $F$ satisfying (2.6) and (2.7). Combining (4.7) and (4.9) gives

$$\lim_{\alpha \downarrow 0} \frac{1}{\sqrt{\alpha}} |g_F(1,\alpha) - \mu_F - 2\sigma_F \sqrt{\alpha}| = 0,$$

as required for (2.8). □

## 5. Growth models.

5.1. *Definitions and statement of shape theorem.* Recall that first-passage quantities $S$ and $h$, analogous to $T$ and $g$, were defined at (2.9) and (2.10).

Define $B(t)$, the *last-passage shape at time $t$*, by

$$B(t) = \{\mathbf{x} \in \mathbb{R}_+^d : T(\lfloor \mathbf{x} \rfloor) \leq t\},$$

and define $C(t)$, the *first-passage shape at time* $t$, by
$$C(t) = \{\mathbf{x} \in \mathbb{R}_+^d : S(\lfloor \mathbf{x} \rfloor) \leq t\}.$$

Both $B(t)$ and $C(t)$ are increasing in $t$, in the sense that for $0 \leq t_1 \leq t_2$, $B(t_1) \subseteq B(t_2)$ and $C(t_1) \subseteq C(t_2)$.

We further define subsets $B$ and $C$ of $\mathbb{R}_+^d$ by
$$B = \{\mathbf{x} : g(\mathbf{x}) \leq 1\},$$
$$C = \{\mathbf{x} : h(\mathbf{x}) \leq 1\}.$$

$B$ is concave [by Proposition 2.1(iv)]; similarly, $C$ is convex. $B$ and $C$ are *asymptotic shapes* for the processes $\{B(t)\}$ and $\{C(t)\}$ in the sense of the following theorem, which is analogous to well-known results for undirected first-passage percolation models (see, e.g., [4, 14]). We give the proof in Section 5.3.

THEOREM 5.1. *Suppose that the weight distribution $F$ satisfies*

(5.1) $$\int_{-\infty}^{0} F(s)^{1/d} ds < \infty$$

*and*

(5.2) $$\int_{0}^{\infty} (1 - F(s))^{1/d} ds < \infty.$$

(i) *Last-passage shape theorem.*
*If $\mathbb{E}_F X > 0$, then for any $\varepsilon > 0$,*

(5.3) $$(1 - \varepsilon)B \subseteq \frac{B(t)}{t} \subseteq (1 + \varepsilon)B$$

*for all sufficiently large $t$, with probability* 1.

(ii) *First-passage shape theorem. If*
$$\inf_{\mathbf{x} \in \mathbb{R}_+^d \setminus \{\mathbf{0}\}} \frac{h(\mathbf{x})}{\|\mathbf{x}\|} > 0,$$

*then for any $\varepsilon > 0$,*
$$(1 - \varepsilon)C \subseteq \frac{C(t)}{t} \subseteq (1 + \varepsilon)C$$

*for all sufficiently large $t$, with probability* 1.

REMARK 5.1. (i) Note that by a subadditivity property for $h$, analogous to the superadditivity property for $g$ in Proposition 2.1(iv), we have that for all $\mathbf{x} \in \mathbb{R}_+^d$, $h(\mathbf{x}) \geq \frac{\|\mathbf{x}\|}{d} h(1, 1, \ldots, 1)$; thus the condition in part (ii) is equivalent to the condition that $h(1, 1, \ldots, 1) > 0$. If the weights are non-negative, then this is implied by the condition that $F(0) < p_c^{(d)}$, where $p_c^{(d)}$ is the critical value for directed percolation in $d$ dimensions; this follows, for example, from the same arguments as the property, noted by Kesten and Su [17], that (in their case for undirected percolation) the critical points for percolation and for "1-percolation" coincide.

(ii) In fact, the moment conditions in Theorem 5.1 are stronger than necessary; for the last-passage case one can replace (5.1) by the condition that $\mathbb{E}|X_-| < \infty$, and for the first-passage case one can replace (5.2) by the condition that $\mathbb{E} X < \infty$.

In particular, combining with the previous remark, for the first-passage model with nonnegative weights it suffices for the limiting shape result that $\mathbb{E} X < \infty$ and $F(0) < p_c^{(d)}$.

To prove the theorem under these weaker conditions, one can follow an approach similar to that used by Cox and Durrett [4] for the undirected first-passage case, making use of the fact that $\mathbb{E} X < \infty \Rightarrow \mathbb{E} \min(X_1, \ldots, X_d)^d < \infty$, where $X_1, \ldots, X_d$ are i.i.d. copies of $X$. However, the proof needs a rather lengthy enumeration of cases and a description of a variety of different sets of "alternative paths"; in addition, the fundamental ideas are already in [4], so we do not include it here.

Under the stronger conditions in Theorem 5.1 (already almost optimal for the last-passage case with nonnegative weights), bounds of the sort given in Lemma 3.5 are available, and the proof is simpler.

5.2. *Discussion of the growth models.* In this section we describe and illustrate various possible types of behavior for the first-passage and last-passage growth processes $B(t)$ and $C(t)$.

We will concentrate on the case where the weights $\{X(\mathbf{z})\}$ are nonnegative. Then for all $t \geq 0$, $B(t)$ and $C(t)$ are connected; in addition, $B(t)$ is a decreasing subset of $\mathbb{R}_+^d$, although $C(t)$ does not generally have a similarly simple property.

If the weight distribution $F$ is exponential (resp. geometric), then the processes $\{B(t), t \geq 0\}$ and $\{C(t), t \geq 0\}$ are Markov in continuous (resp. discrete) time; in the two-dimensional last-passage case this yields the growth model considered in, for example, Rost [25] and Johansson [13]. Simulations of $B(t)$ and $C(t)$ in two dimensions with exponentially distributed weights are given in Figures 1 and 2, and the three-dimensional last-passage case is shown in Figure 7.

In fact $B(t)$ and $C(t)$ are also Markov in discrete time when the weights are Bernoulli (taking values 0 and 1).

We now discuss various ways in which the shape result in Theorem 5.1 may fail.

First, the last-passage case. Note that if $g(\mathbf{x}) = \infty$ for some $\mathbf{x}$ in the interior of $\mathbb{R}_+^d$, then (by a simple superadditivity argument) $g = \infty$ throughout the interior. A sufficient condition for this to occur is that $\mathbb{E} X^d = \infty$. Then the growth of $B(t)$ in any interior direction is sublinear in $t$; on a linear scale, the asymptotic shape collapses into the boundary (or even to the origin alone if $\mathbb{E} X = \infty$). An example is illustrated in Figure 3 for a distribution with finite mean but infinite variance, with $d = 2$.

The undirected first-passage case on $\mathbb{Z}^d$ is comparable. Here again there are just two possibilities: either the shape function is 0 everywhere (in which case the asymptotic shape is essentially the whole of $\mathbb{R}^d$), or the shape function is nonzero everywhere.

For the directed first-passage case, in contrast, the different behaviors can co-exist. If $p_c^{(d)} < \mathbb{P}(X = 0) < 1$, then $h = 0$ for some cone around the direction $\mathbf{1}$ (this is the cone in which "oriented percolation of sites with weight 0" occurs), but $h > 0$ elsewhere. In this case there will a.s. be some infinite path starting at the origin which has finite total weight, and the set $C(t)$ will have infinite size at some finite time. The shape $C = \{\mathbf{x} : h(\mathbf{x}) \leq 1\}$ is noncompact, but not equal to the whole of $\mathbb{R}_+^d$, and the shape theorem does not apply as given. See Figure 4 for a simulation of such a case.

Note that if again $X$ attains its minimum value with probability more than $p_c^{(d)}$, but this minimum is now greater than 0, then $h(\mathbf{x})/\|\mathbf{x}\|$ is constant on some cone around the direction $\mathbf{1}$, and the boundary of $C$ has a flat section— see Figure 5. This is the same phenomenon observed by Durrett and Liggett [8] for the undirected case.

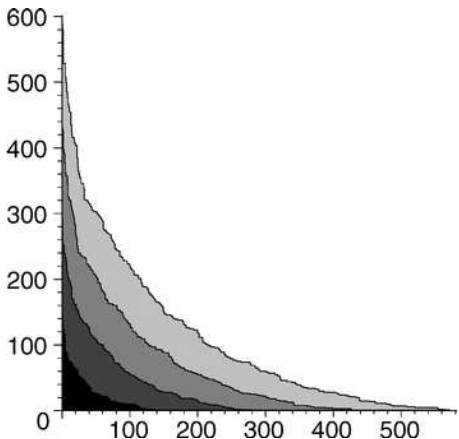 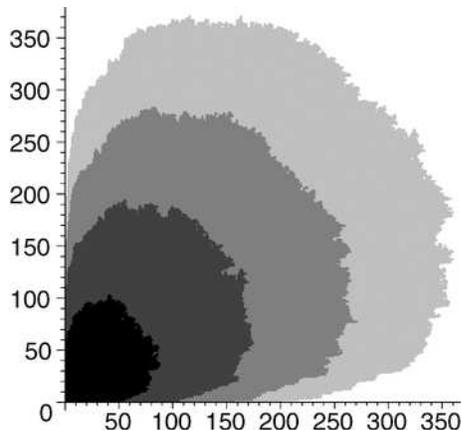

Fig. 1. *Simulation of the last-passage process for $d = 2$ and $F$ exponential with mean 1. The sets $B(t)$ are shown for $t = 150$ (darkest), $300, 450, 600$ (lightest). Here the asymptotic shape $B$ is known to be $\{(x, y) \in \mathbb{R}_+^2 : \sqrt{x} + \sqrt{y} \leq 1\}$.*

Fig. 2. *First-passage process $C(t)$, $t = 50, 100, 150, 200$ for $d = 2$ and $F$ exponential with mean 1.*

If $\mathbb{E} X = \infty$, then $h(\mathbf{e}_1) = \infty$, and the shape theorem cannot hold as stated, but an amended version for a cone excluding the boundaries may hold. However, if (only slightly more strongly) $\mathbb{E} \min(X_1, \ldots, X_d)^d = \infty$, then the result fails more fundamentally. Still the limit $h(\mathbf{x}) = \lim_{n \to \infty} n^{-1} S(\lfloor n\mathbf{x} \rfloor)$ exists and is finite and constant a.s. for any $\mathbf{x}$ in the interior of $\mathbb{R}_+^d$, and one can define the asymptotic shape $C$ as before. However, it is no longer the case that the convergence is a.s. uniform on compact subsets of $\mathbb{R}_+^d$; in effect, the "holes" in the shape $C(t)$ persist, as seen in Figure 6. The same sort of behavior would occur for the undirected first-passage model when $\mathbb{E} \min(X_1, \ldots, X_{2d})^d = \infty$; see, for example, related discussions in [4].

5.3. *Proof of Theorem* 5.1. Note first that for part (i) of Theorem 5.1, the condition $\mathbb{E} X > 0$ implies that $\inf_{\mathbf{x} \in \mathbb{R}_+^d \setminus \{\mathbf{0}\}} \frac{g(\mathbf{x})}{\|\mathbf{x}\|} > 0$ [since, by superadditivity, $g(\mathbf{x}) \geq \|\mathbf{x}\| g(1, 0, \ldots, 0) = \|\mathbf{x}\| \mathbb{E} X$].

Note also that, by replacing the weights $\{X(\mathbf{z})\}$ by $\{-X(\mathbf{z})\}$, part (ii) can be rewritten as follows in terms of last-passage rather than first-passage quantities:

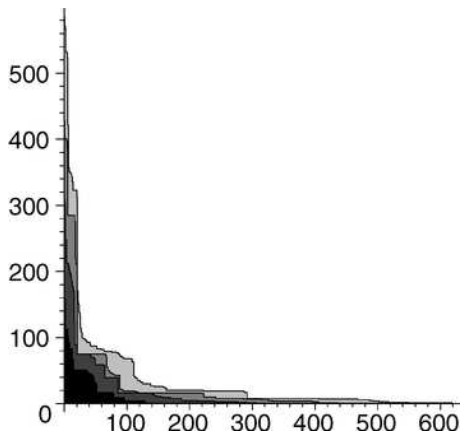

FIG. 3. $B(t)$, $t = 150, 300, 450, 600$, with $d = 2$ and the Pareto distribution $F(x) = \min(0, 1 - (3x)^{-3/2})$, which has mean 1 but infinite variance. The asymptotic shape $B$ consists only of two lines, between the origin and $(0, 1)$ and between the origin and $(1, 0)$.

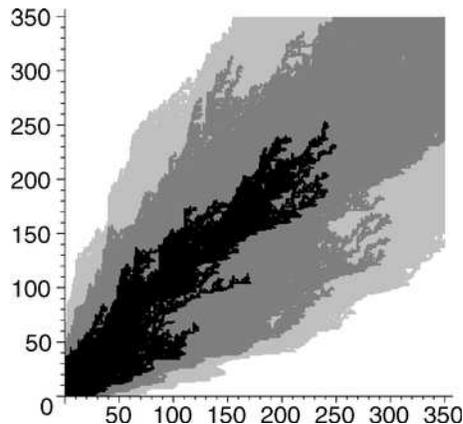

FIG. 4. $C(t)$, $t = 18, 36, 54$, with $d = 2$ and $X = 0$ w.p. 0.645 and $X = (0.355)^{-1}$ w.p. 0.355 (showing only the intersection with the box $[0, 350]^2$). Here $h(\mathbf{x}) = 0$ for some (though not all) $\mathbf{x}$; thus the asymptotic shape is noncompact (but not the whole of $\mathbb{R}_+^d$), and a.s. the set $C(t)$ will be infinite for some $t$.

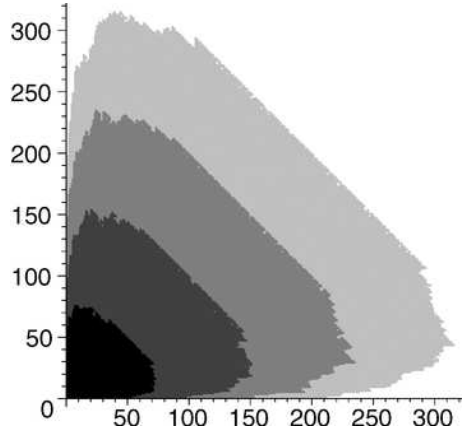

Fig. 5. $C(t)$, $t = 50, 100, 150, 200$, with $d = 2$ and $X = 0.5$ w.p. 0.8, $X = 3$ w.p. 0.2.

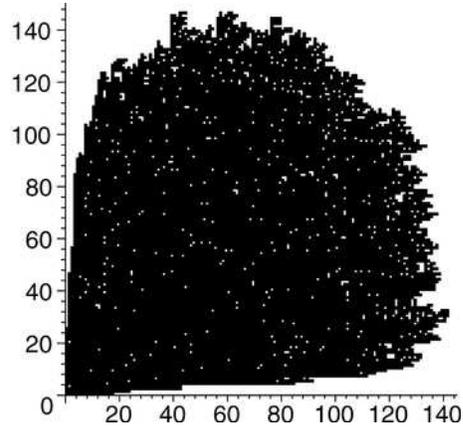

Fig. 6. $C(400)$ for $d = 2$ and the Pareto distribution $F(x) = \min(0, 1 - x^{-3/4})$, which has infinite mean.

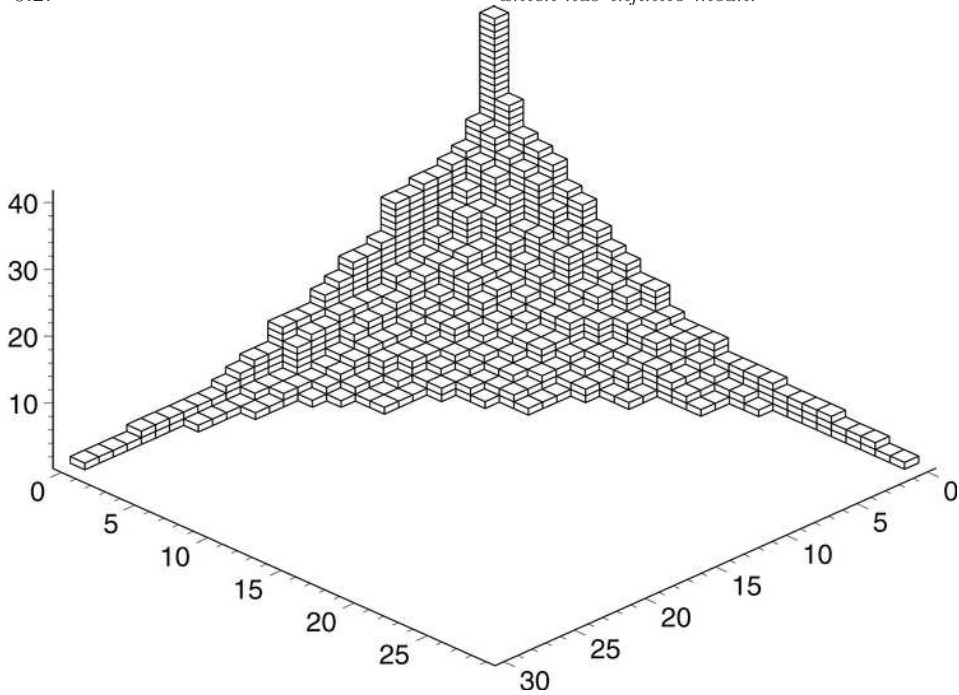

Fig. 7. $B(30)$ for $d = 3$ and $F$ exponential with mean 1.

If

$$\sup_{\mathbf{x} \in \mathbb{R}_+^d \setminus \{\mathbf{0}\}} \frac{g(\mathbf{x})}{\|\mathbf{x}\|} < 0,$$

*then for any $\varepsilon > 0$,*

(5.4) $$(1-\varepsilon)B^- \subseteq \frac{B^-(t)}{t} \subseteq (1+\varepsilon)B^-$$

*for all sufficiently large $t$, with probability 1, where*

$$B^-(t) = \{\mathbf{x} \in \mathbb{R}_+^d : T(\lfloor \mathbf{x} \rfloor) \geq -t\}$$

*and*

$$B^- = \{\mathbf{x} : g(\mathbf{x}) \geq -1\}.$$

By Theorem 2.3, we know that $|g(\mathbf{z})| < \infty$ for all $\mathbf{z}$; then (5.3) and (5.4) are immediately implied by the following property: for any $\varepsilon > 0$, there are a.s. only finitely many $\mathbf{z} \in \mathbb{Z}_+^d$ such that $|T(\mathbf{z}) - g(\mathbf{z})| \geq \varepsilon |g(\mathbf{z})|$.

Since we assume in both cases that $|g(\mathbf{z})|/\|\mathbf{z}\|$ is bounded away from 0, we have $|g(\mathbf{z})| \geq \|\mathbf{z}\| \inf_{\mathbf{z}'} |g(\mathbf{z}')|/\|\mathbf{z}'\|$, so in fact it is enough to show that for any $\varepsilon > 0$, there are a.s. only finitely many $\mathbf{z}$ such that $|T(\mathbf{z}) - g(\mathbf{z})| \geq \varepsilon \|\mathbf{z}\|$.

This follows immediately from Lemmas 5.2–5.5, which we state immediately and then prove in turn.

LEMMA 5.2. *Suppose $F$ satisfies $\int_{-\infty}^{0} F(s)^{1/d} ds < \infty$ and $\int_0^\infty (1-F(s))^{1/d} ds < \infty$, and let $\varepsilon > 0$. If $L$ is large enough, then with probability 1,*

$$|T(\mathbf{z}) - T^{(L)}(\mathbf{z})| < \varepsilon \|\mathbf{z}\|$$

*for all except finitely many $\mathbf{z} \in \mathbb{Z}_+^d$.*

LEMMA 5.3. *Let $\varepsilon > 0$ and $L > 0$. With probability 1,*

$$|T^{(L)}(\mathbf{z}) - \mathbb{E}\,T^{(L)}(\mathbf{z})| < \varepsilon \|\mathbf{z}\|,$$

*for all except finitely many $\mathbf{z} \in \mathbb{Z}_+^d$.*

LEMMA 5.4. *Let $\varepsilon > 0$ and $L > 0$. Then for all except finitely many $\mathbf{z} \in \mathbb{Z}_+^d$,*

$$|\mathbb{E}\,T^{(L)}(\mathbf{z}) - g^{(L)}(\mathbf{z})| < \varepsilon \|\mathbf{z}\|.$$

LEMMA 5.5. *Suppose $F$ satisfies $\int_{-\infty}^{0} F(s)^{1/d} ds < \infty$ and $\int_0^\infty (1-F(s))^{1/d} ds < \infty$, and let $\varepsilon > 0$. If $L$ is sufficiently large, then for all $\mathbf{z} \in \mathbb{Z}_+^d$,*

$$|g^{(L)}(\mathbf{z}) - g(\mathbf{z})| < \varepsilon \|\mathbf{z}\|.$$

PROOF OF LEMMA 5.2. Let $L$ be large enough that $c \int_L^\infty (1 - F(s))^{1/d} ds < \varepsilon/2$ and $c \int_{-\infty}^{-L} F(s)^{1/d} ds < \varepsilon/2$, where $c$ is the constant in Lemma 3.5.

Let $\mathbf{z} \in \mathbb{Z}_+^d$; for some $\pi^* \in \Pi[\mathbf{z}]$, we have

$$T(\mathbf{z}) - T^{(L)}(\mathbf{z}) = \sum_{\mathbf{v} \in \pi^*} [X(\mathbf{v}) - X^{(L)}(\mathbf{v})],$$

so that

(5.5)
$$|T(\mathbf{z}) - T^{(L)}(\mathbf{z})| \leq \sum_{\mathbf{v} \in \pi^*} [X(\mathbf{v}) - L]_+ + \sum_{\mathbf{v} \in \pi^*} [-L - X(\mathbf{v})]_+$$
$$\leq \max_{\pi \in \Pi[\mathbf{z}]} \sum_{\mathbf{v} \in \pi} V^{(L)}(\mathbf{v}) + \max_{\pi \in \Pi[\mathbf{z}]} \sum_{\mathbf{v} \in \pi} W^{(L)}(\mathbf{v}),$$

where we define $V^{(L)}(\mathbf{v}) = [X(\mathbf{v}) - L]_+$ and $W^{(L)}(\mathbf{v}) = [-L - X(\mathbf{v})]_+$.

Note that $\{V^{(L)}(\mathbf{v}), \mathbf{v} \in \mathbb{Z}_+^d\}$ are i.i.d. with common distribution $F_V^{(L)}$, where $F_V^{(L)}(x) = 0$ for $x < 0$ and $F_V^{(L)}(x) = F(L + x)$ for $x \geq 0$.

Similarly, $\{W^{(L)}(\mathbf{v}), \mathbf{v} \in \mathbb{Z}_+^d\}$ are i.i.d. with common distribution $F_W^{(L)}$, where $F_W^{(L)}(x) = 0$ for $x < 0$ and $F_W^{(L)}(x) = 1 - F(-L - x)$ for $x \geq 0$.

From Lemma 3.5 (applied to $\{V^{(L)}(\mathbf{v})\}$ and $F_V^{(L)}$ rather than to $\{X(\mathbf{v})\}$ and $F$), we then have that, with probability 1,

$$\limsup_{n \to \infty} \frac{1}{n} \max_{\mathbf{z}: \|\mathbf{z}\| \leq n} \max_{\pi \in \Pi[\mathbf{z}]} \sum_{\mathbf{v} \in \pi} V^{(L)}(\mathbf{v}) \leq c \int_0^\infty (1 - F_V^{(L)}(s))^{1/d} ds$$
$$= c \int_L^\infty (1 - F(s))^{1/d} ds$$
$$< \varepsilon/2.$$

In particular, there are a.s. only finitely many $\mathbf{z} \in \mathbb{Z}_+^d$ such that

$$\max_{\pi \in \Pi[\mathbf{z}]} \sum_{\mathbf{v} \in \pi} V^{(L)}(\mathbf{v}) \geq \frac{\varepsilon}{2} \|\mathbf{z}\|.$$

Applying Lemma 3.5 to $\{W^{(L)}(\mathbf{v})\}$ and $F_W^{(L)}$ in the same way, one obtains that there are a.s. only finitely many $\mathbf{z} \in \mathbb{Z}_+^d$ such that

$$\max_{\pi \in \Pi[\mathbf{z}]} \sum_{\mathbf{v} \in \pi} W^{(L)}(\mathbf{v}) \geq \frac{\varepsilon}{2} \|\mathbf{z}\|.$$

Thus from (5.5), there are a.s. only finitely many $\mathbf{z} \in \mathbb{Z}_+^d$ such that

$$|T(\mathbf{z}) - T^{(L)}(\mathbf{z})| \geq \varepsilon \|\mathbf{z}\|,$$

as required. □

PROOF OF LEMMA 5.3. All the paths in $\Pi[\mathbf{z}]$ have length $\|\mathbf{z}\|$, and the weights $X^{(L)}(\mathbf{v})$ have absolute value no greater than $L$. Hence we may apply the concentration inequality in Lemma 3.1 to give

$$\mathbb{P}(|T^{(L)}(\mathbf{z}) - \mathbb{E}\, T^{(L)}(\mathbf{z})| \geq \varepsilon \|\mathbf{z}\|) \leq \exp\left(-\frac{(\varepsilon\|\mathbf{z}\|)^2}{64\|\mathbf{z}\|L^2} + 64\right)$$

$$= \exp\left(-\frac{\varepsilon^2\|\mathbf{z}\|}{64L^2} + 64\right).$$

For $n \in \mathbb{Z}_+$, there are certainly no more than $(n+1)^d$ points $\mathbf{z}$ such that $\|\mathbf{z}\| = n$; thus

$$\sum_{\mathbf{z} \in \mathbb{Z}_+^d} \mathbb{P}(|T^{(L)}(\mathbf{z}) - \mathbb{E}\, T^{(L)}(\mathbf{z})| \geq \varepsilon\|\mathbf{z}\|) \leq \sum_{n \in \mathbb{Z}_+^d} (n+1)^d \exp\left(-\frac{\varepsilon^2 n}{64L^2} + 64\right)$$

$$< \infty,$$

and Borel–Cantelli yields the desired result. $\square$

PROOF OF LEMMA 5.4. From the definition of $g$ and so of $g^{(L)}$, we have $\mathbb{E}\, T^{(L)}(\mathbf{z}) \leq g^{(L)}(\mathbf{z})$ for all $\mathbf{z}$, so we need to show that $\mathbb{E}\, T^{(L)}(\mathbf{z}) > g^{(L)}(\mathbf{z}) - \varepsilon\|\mathbf{z}\|$, except for finitely many $\mathbf{z}$.

Fix $a > 0$. The distribution $F^{(L)}$ has bounded support, so certainly satisfies (2.4) and (2.5); thus by Theorem 2.3, $g^{(L)}$ is continuous on $\mathbb{R}_+^d$, and hence is uniformly continuous on the compact subset $\{\mathbf{x} \in \mathbb{R}_+^d : \|\mathbf{x}\| \leq 2d\}$.

So choose $u < \min(1, a)$ such that whenever $\|\mathbf{x}\| \leq d$ and $\|\mathbf{x} - \mathbf{x}'\| \leq ud$, then $|g^{(L)}(\mathbf{x}) - g^{(L)}(\mathbf{x}')| \leq a$.

Now let $C$ be the set

$$\left\{u\mathbf{r}, \mathbf{r} \in \left\{0, 1, \ldots, \left\lfloor \frac{1}{u} \right\rfloor\right\}^d\right\}.$$

$C$ is a finite subset of $\mathbb{R}_+^d$, and for every $\mathbf{y} \in C$, we have [by Proposition 2.1(i)],

$$\frac{\mathbb{E}\, T^{(L)}(\lfloor n\mathbf{y} \rfloor)}{n} \to g^{(L)}(\mathbf{y}) \qquad \text{as } n \to \infty.$$

Hence there is $N = N(a)$ such that, for all $n \geq N$ and all $\mathbf{y} \in C$,

$$\mathbb{E}\, T^{(L)}(\lfloor n\mathbf{y} \rfloor) \geq n(g^{(L)}(\mathbf{y}) - a).$$

Let $\mathbf{z}$ satisfy $\max z_i \geq N$. Define

$$\mathbf{y} = u\left\lfloor \frac{1}{u}\frac{\mathbf{z}}{\max z_i} \right\rfloor;$$

then $\mathbf{y} \in C$, with $(\max z_i)\mathbf{y} \leq \mathbf{z}$, with $\|\mathbf{y}\| \leq d$ and with

$$\left\| \frac{\mathbf{z}}{\max z_i} - \mathbf{y} \right\| \leq ud \leq ad.$$

Using first superadditivity, then the fact that all the weights $\{X^{(L)}(\mathbf{z})\}$ are no smaller than $-L$, then the continuity bounds above, we obtain

$$\begin{aligned}
\mathbb{E}\,T^{(L)}(\mathbf{z}) &\geq \mathbb{E}\,T^{(L)}(\lfloor (\max z_i)\mathbf{y}\rfloor) + \mathbb{E}\,T^{(L)}(\mathbf{z} - \lfloor (\max z_i)\mathbf{y}\rfloor) \\
&\geq \mathbb{E}\,T^{(L)}(\lfloor (\max Z_i)\mathbf{y}\rfloor) - L\|\mathbf{z} - \lfloor (\max Z_i)\mathbf{y}\rfloor\| \\
&\geq (\max z_i)(g^{(L)}(\mathbf{y}) - a) - L(\|\mathbf{z} - (\max z_i)\mathbf{y}\| + d) \\
&= g^{(L)}(\mathbf{z}) - (\max z_i)\left\{\left[g^{(L)}\!\left(\frac{\mathbf{z}}{\max z_i}\right) - g^{(L)}(\mathbf{y})\right] \right. \\
&\qquad\qquad\qquad\qquad \left. + a + L\left\|\frac{\mathbf{z}}{\max z_i} - \mathbf{y}\right\| + \frac{Ld}{\max z_i}\right\} \\
&\geq g^{(L)}(\mathbf{z}) - (\max z_i)\left\{a + a + Lad + \frac{Ld}{\max z_i}\right\}.
\end{aligned}$$

Hence if $a < 2\varepsilon^{-1}(2+Ld)$, then for all $\mathbf{z}$ with $\max z_i \geq \max(N(a), 2Ld/\varepsilon)$, we have

$$\begin{aligned}
\mathbb{E}\,T^{(L)}(\mathbf{z}) &> g^{(L)}(\mathbf{z}) - (\max z_i)\varepsilon \\
&\geq g^{(L)}(\mathbf{z}) - \varepsilon\|\mathbf{z}\|,
\end{aligned}$$

as required. $\square$

PROOF OF LEMMA 5.5. Under the moment conditions on $F$, the result follows immediately from Lemma 3.6.

This completes the proof of Theorem 5.1. $\square$

LIAFA
CNRS et Université Paris 7
case 7014
2 place Jussieu
75251 Paris Cedex 05
France
e-mail: James.Martin@liafa.jussieu.fr